\documentclass{gtart}
\def\ifplaintex{\expandafter\ifx\csname documentclass\endcsname\relax}

\def\gtp{{\mathsurround=0pt\it $\cal G\mskip-2mu$eometry \&\ 
$\cal T\!\!$opology $\cal P\!$ublications}}  % GT publications

\def\recd{{\small Received:\qua\receiveddate\ifx\reviseddate\relax
\else\qquad Revised:\qua\reviseddate\fi\par}} 

\def\lognumber#1{\def\thelognumber{#1}}
\def\volumenumber#1{\def\thevolumenumber{#1}}
\def\volumeyear#1{\def\thevolumeyear{#1}}
\def\papernumber#1{\def\thepapernumber{#1}}
\def\pagenumbers#1#2{\def\startpage{#1}\def\finishpage{#2}}
\def\published#1{\def\publishdate{#1}}

\def\received#1{\def\receiveddate{#1}}
\def\revised#1{\def\reviseddate{#1}}
\def\accepted#1{\def\accepteddate{#1}}

\long\def\asciiabstract#1{\long\def\theasciiabstract{#1}}
\def\asciikeywords#1{\def\theasciikeywords{#1}}

\let\\\par\let\thelognumber\relax\let\thevolumenumber\relax
\let\thepapernumber\relax\let\thevolumeyear\relax\let\startpage\relax
\let\finishpage\relax\let\publishdate\relax\let\receiveddate\relax
\let\reviseddate\relax\let\accepteddate\relax\let\theasciititle\relax
\let\theasciiauthors\relax
\let\theasciiabstract\relax\let\theasciikeywords\relax

\let\theasciiemail\relax

\ifplaintex
\font\logobig=cmssbx10 scaled 3836
\font\logomed=cmssbx10 scaled 2557
\else
\font\logobig=cmssbx10 scaled 4200
\font\logomed=cmssbx10 scaled 2800
\fi

\long\def\makeagttitle{   %%% start of definition of \makeagttitle
\count0=\startpage
\agt\hfill      %   Journal title (top left) 
\hbox to 45truept{\vbox to 0pt{\vglue -13truept{\logomed A\kern -.37em{\logobig 
T}\kern -.38em G}\vss}\hss}
\break
{\small Volume \thevolumenumber\ (\thevolumeyear)
\startpage--\finishpage\nl
Published: \publishdate}

\vglue .25truein

{\parskip=0pt\leftskip 0pt plus
1fil\def\\{\par\smallskip}{\Large\bf\thetitle}\par\medskip} \vglue
0.05truein

{\parskip=0pt\leftskip 0pt plus 1fil\def\\{\par}{\sc\theauthors}
\par\medskip}%
 
\vglue 0.03truein 

{\small\leftskip 25truept\rightskip 25truept{\bf Abstract}\stdspace\theabstract

{\bf AMS Classification}\stdspace\theprimaryclass
\ifx\thesecondaryclass\relax\else; \thesecondaryclass\fi\par
{\bf Keywords}\stdspace \thekeywords\par}\vglue 7truept

}   %%%% end of definition of \makeagttitle

\ifplaintex
\font\phead=cmsl9 scaled 950
\font\pnum=cmbx10 scaled 913
\font\pfoot=cmsl9 scaled 950
\headline{\vbox to 0pt{\vskip -4.5mm\line{\small\phead\ifnum
\count0=\startpage ISSN 1472-2739 (on-line) 1472-2747 (printed)
\hfill {\pnum\folio}\else\ifodd\count0\def\\{ }% 
\ifx\theshorttitle\relax\thetitle\else\theshorttitle\fi\hfill{\pnum\folio}
\else\def\\{ and }{\pnum\folio}\hfill\ifx\theshortauthors\relax\theauthors
\else\theshortauthors\fi\fi\fi}\vss}}
\footline{\vbox to 0pt{\vglue 0mm\line{\small\pfoot\ifnum\count0=\startpage
\copyright\ \gtp\hfill\else
\agt, Volume \thevolumenumber\ (\thevolumeyear)\hfill\fi}\vss}}
\else
\font=cmsl9 scaled 1050
\font\lnum=cmbx10 
\font=cmsl9 scaled 1050
\makeatletter
\def\@oddhead{{\small\lhead\ifnum\count0=\startpage ISSN 1472-2739 
(on-line) 1472-2747 (printed)\hfill {\lnum\number\count0}\else\ifodd\count0
\def\\{ }\ifx\theshorttitle\relax \thetitle \else\theshorttitle\fi\hfill
{\lnum\number\count0}\else\def\\{ and }{\lnum\number\count0}
\hfill\ifx\theshortauthors\relax 
\theauthors\else\theshortauthors\fi\fi\fi}}\def\@evenhead{\@oddhead}
\def\@oddfoot{\small\lfoot\ifnum\count0=\startpage\copyright\ \gtp\hfill\else
\agt, Volume \thevolumenumber\ (\thevolumeyear)\hfill\fi}
\def\@evenfoot{\@oddfoot}
\makeatother
\fi
\let\maketitlepage\makeagttitle
\let\makeshorttitle\maketitlepage
\let\maketitle\maketitlepage

\newwrite\gtoutfile
\long\gdef\makeheadfile{  %%% start of definition of \makeheadfile
{\def\\{, }\def\s{ }
\immediate\openout\gtoutfile head.xxx
\immediate\write\gtoutfile{To: math@arxiv.org}
\immediate\write\gtoutfile{Subject: put OR rep NNNNN:ppppp}
\immediate\write\gtoutfile{--text follows this line--}
\immediate\write\gtoutfile{Proxy-for: \ifx\theasciiauthors\relax
\theauthors\else\theasciiauthors\fi\s<\ifx\theasciiemail\relax\theemail\else\theasciiemail\fi>}
\immediate\write\gtoutfile{\noexpand\\}
\immediate\write\gtoutfile{Authors: \ifx\theasciiauthors\relax
\theauthors\else\theasciiauthors\fi}
{\def\\{ }\immediate\write\gtoutfile{Title: \ifx\theasciititle\relax
\thetitle\else\theasciititle\fi}}
\immediate\write\gtoutfile{Subj-class: GT or SG, GR etc}
\immediate\write\gtoutfile{MSC-class: \theprimaryclass\ifx\thesecondaryclass\relax\else, \thesecondaryclass\fi}
\immediate\write\gtoutfile{Journal-ref: Algebr. Geom. Topol. \thevolumenumber\s
(\thevolumeyear) \startpage-\finishpage}
\immediate\write\gtoutfile{Comments: Published by Algebraic and
Geometric Topology at}
\immediate\write\gtoutfile{\s\s\s  http://www.maths.warwick.ac.uk/agt/AGTVol\thevolumenumber/agt-\thevolumenumber-\thepapernumber.abs.html}
\immediate\write\gtoutfile{\noexpand\\}
\immediate\write\gtoutfile{}
\ifx\theasciiabstract\relax
\immediate\write\gtoutfile{\theabstract}\else
\immediate\write\gtoutfile{\theasciiabstract}\fi
\immediate\write\gtoutfile{}
\immediate\write\gtoutfile{\noexpand\\}
\immediate\write\gtoutfile{}
\immediate\closeout\gtoutfile}}  %%% end of definition of \makeheadfile

\def\maketitlepage{\makeagttitle\makeheadfile}
\let\makeshorttitle\maketitlepage
\let\maketitle\maketitlepage

\def\ifplaintex{\expandafter\ifx\csname documentclass\endcsname\relax}

\def\gtp{{\mathsurround=0pt\it $\cal G\mskip-2mu$eometry \&\ 
$\cal T\!\!$opology $\cal P\!$ublications}}  % GT publications

\def\recd{{\small Received:\qua\receiveddate\ifx\reviseddate\relax
\else\qquad Revised:\qua\reviseddate\fi\par}} 

\def\lognumber#1{\def\thelognumber{#1}}
\def\volumenumber#1{\def\thevolumenumber{#1}}
\def\volumeyear#1{\def\thevolumeyear{#1}}
\def\papernumber#1{\def\thepapernumber{#1}}
\def\pagenumbers#1#2{\def\startpage{#1}\def\finishpage{#2}}
\def\published#1{\def\publishdate{#1}}

\def\received#1{\def\receiveddate{#1}}
\def\revised#1{\def\reviseddate{#1}}
\def\accepted#1{\def\accepteddate{#1}}

\long\def\asciiabstract#1{\long\def\theasciiabstract{#1}}
\def\asciikeywords#1{\def\theasciikeywords{#1}}

\let\\\par\let\thelognumber\relax\let\thevolumenumber\relax
\let\thepapernumber\relax\let\thevolumeyear\relax\let\startpage\relax
\let\finishpage\relax\let\publishdate\relax\let\receiveddate\relax
\let\reviseddate\relax\let\accepteddate\relax\let\theasciititle\relax
\let\theasciiauthors\relax
\let\theasciiabstract\relax\let\theasciikeywords\relax

\let\theasciiemail\relax

\ifplaintex
\font\logobig=cmssbx10 scaled 3836
\font\logomed=cmssbx10 scaled 2557
\else
\font\logobig=cmssbx10 scaled 4200
\font\logomed=cmssbx10 scaled 2800
\fi

\long\def\makeagttitle{   %%% start of definition of \makeagttitle
\count0=\startpage
\agt\hfill      %   Journal title (top left) 
\hbox to 45truept{\vbox to 0pt{\vglue -13truept{\logomed A\kern -.37em{\logobig 
T}\kern -.38em G}\vss}\hss}
\break
{\small Volume \thevolumenumber\ (\thevolumeyear)
\startpage--\finishpage\nl
Published: \publishdate}

\vglue .25truein

{\parskip=0pt\leftskip 0pt plus
1fil\def\\{\par\smallskip}{\Large\bf\thetitle}\par\medskip} \vglue
0.05truein

{\parskip=0pt\leftskip 0pt plus 1fil\def\\{\par}{\sc\theauthors}
\par\medskip}%
 
\vglue 0.03truein 

{\small\leftskip 25truept\rightskip 25truept{\bf Abstract}\stdspace\theabstract

{\bf AMS Classification}\stdspace\theprimaryclass
\ifx\thesecondaryclass\relax\else; \thesecondaryclass\fi\par
{\bf Keywords}\stdspace \thekeywords\par}\vglue 7truept

}   %%%% end of definition of \makeagttitle

\ifplaintex
\font\phead=cmsl9 scaled 950
\font\pnum=cmbx10 scaled 913
\font\pfoot=cmsl9 scaled 950
\headline{\vbox to 0pt{\vskip -4.5mm\line{\small\phead\ifnum
\count0=\startpage ISSN 1472-2739 (on-line) 1472-2747 (printed)
\hfill {\pnum\folio}\else\ifodd\count0\def\\{ }% 
\ifx\theshorttitle\relax\thetitle\else\theshorttitle\fi\hfill{\pnum\folio}
\else\def\\{ and }{\pnum\folio}\hfill\ifx\theshortauthors\relax\theauthors
\else\theshortauthors\fi\fi\fi}\vss}}
\footline{\vbox to 0pt{\vglue 0mm\line{\small\pfoot\ifnum\count0=\startpage
\copyright\ \gtp\hfill\else
\agt, Volume \thevolumenumber\ (\thevolumeyear)\hfill\fi}\vss}}
\else
\font=cmsl9 scaled 1050
\font\lnum=cmbx10 
\font=cmsl9 scaled 1050
\makeatletter
\def\@oddhead{{\small\lhead\ifnum\count0=\startpage ISSN 1472-2739 
(on-line) 1472-2747 (printed)\hfill {\lnum\number\count0}\else\ifodd\count0
\def\\{ }\ifx\theshorttitle\relax \thetitle \else\theshorttitle\fi\hfill
{\lnum\number\count0}\else\def\\{ and }{\lnum\number\count0}
\hfill\ifx\theshortauthors\relax 
\theauthors\else\theshortauthors\fi\fi\fi}}\def\@evenhead{\@oddhead}
\def\@oddfoot{\small\lfoot\ifnum\count0=\startpage\copyright\ \gtp\hfill\else
\agt, Volume \thevolumenumber\ (\thevolumeyear)\hfill\fi}
\def\@evenfoot{\@oddfoot}
\makeatother
\fi
\let\maketitlepage\makeagttitle
\let\makeshorttitle\maketitlepage
\let\maketitle\maketitlepage

\newwrite\gtoutfile
\long\gdef\makeheadfile{  %%% start of definition of \makeheadfile
{\def\\{, }\def\s{ }
\immediate\openout\gtoutfile head.xxx
\immediate\write\gtoutfile{To: math@arxiv.org}
\immediate\write\gtoutfile{Subject: put OR rep NNNNN:ppppp}
\immediate\write\gtoutfile{--text follows this line--}
\immediate\write\gtoutfile{Proxy-for: \ifx\theasciiauthors\relax
\theauthors\else\theasciiauthors\fi\s<\ifx\theasciiemail\relax\theemail\else\theasciiemail\fi>}
\immediate\write\gtoutfile{\noexpand\\}
\immediate\write\gtoutfile{Authors: \ifx\theasciiauthors\relax
\theauthors\else\theasciiauthors\fi}
{\def\\{ }\immediate\write\gtoutfile{Title: \ifx\theasciititle\relax
\thetitle\else\theasciititle\fi}}
\immediate\write\gtoutfile{Subj-class: GT or SG, GR etc}
\immediate\write\gtoutfile{MSC-class: \theprimaryclass\ifx\thesecondaryclass\relax\else, \thesecondaryclass\fi}
\immediate\write\gtoutfile{Journal-ref: Algebr. Geom. Topol. \thevolumenumber\s
(\thevolumeyear) \startpage-\finishpage}
\immediate\write\gtoutfile{Comments: Published by Algebraic and
Geometric Topology at}
\immediate\write\gtoutfile{\s\s\s  http://www.maths.warwick.ac.uk/agt/AGTVol\thevolumenumber/agt-\thevolumenumber-\thepapernumber.abs.html}
\immediate\write\gtoutfile{\noexpand\\}
\immediate\write\gtoutfile{}
\ifx\theasciiabstract\relax
\immediate\write\gtoutfile{\theabstract}\else
\immediate\write\gtoutfile{\theasciiabstract}\fi
\immediate\write\gtoutfile{}
\immediate\write\gtoutfile{\noexpand\\}
\immediate\write\gtoutfile{}
\immediate\closeout\gtoutfile}}  %%% end of definition of \makeheadfile

\def\maketitlepage{\makeagttitle\makeheadfile}
\let\makeshorttitle\maketitlepage
\let\maketitle\maketitlepage

\def\ifplaintex{\expandafter\ifx\csname documentclass\endcsname\relax}

\def\gtp{{\mathsurround=0pt\it $\cal G\mskip-2mu$eometry \&\ 
$\cal T\!\!$opology $\cal P\!$ublications}}  % GT publications

\def\recd{{\small Received:\qua\receiveddate\ifx\reviseddate\relax
\else\qquad Revised:\qua\reviseddate\fi\par}} 

\def\lognumber#1{\def\thelognumber{#1}}
\def\volumenumber#1{\def\thevolumenumber{#1}}
\def\volumeyear#1{\def\thevolumeyear{#1}}
\def\papernumber#1{\def\thepapernumber{#1}}
\def\pagenumbers#1#2{\def\startpage{#1}\def\finishpage{#2}}
\def\published#1{\def\publishdate{#1}}

\def\received#1{\def\receiveddate{#1}}
\def\revised#1{\def\reviseddate{#1}}
\def\accepted#1{\def\accepteddate{#1}}

\long\def\asciiabstract#1{\long\def\theasciiabstract{#1}}
\def\asciikeywords#1{\def\theasciikeywords{#1}}

\let\\\par\let\thelognumber\relax\let\thevolumenumber\relax
\let\thepapernumber\relax\let\thevolumeyear\relax\let\startpage\relax
\let\finishpage\relax\let\publishdate\relax\let\receiveddate\relax
\let\reviseddate\relax\let\accepteddate\relax\let\theasciititle\relax
\let\theasciiauthors\relax
\let\theasciiabstract\relax\let\theasciikeywords\relax

\let\theasciiemail\relax

\ifplaintex
\font\logobig=cmssbx10 scaled 3836
\font\logomed=cmssbx10 scaled 2557
\else
\font\logobig=cmssbx10 scaled 4200
\font\logomed=cmssbx10 scaled 2800
\fi

\long\def\makeagttitle{   %%% start of definition of \makeagttitle
\count0=\startpage
\agt\hfill      %   Journal title (top left) 
\hbox to 45truept{\vbox to 0pt{\vglue -13truept{\logomed A\kern -.37em{\logobig 
T}\kern -.38em G}\vss}\hss}
\break
{\small Volume \thevolumenumber\ (\thevolumeyear)
\startpage--\finishpage\nl
Published: \publishdate}

\vglue .25truein

{\parskip=0pt\leftskip 0pt plus
1fil\def\\{\par\smallskip}{\Large\bf\thetitle}\par\medskip} \vglue
0.05truein

{\parskip=0pt\leftskip 0pt plus 1fil\def\\{\par}{\sc\theauthors}
\par\medskip}%
 
\vglue 0.03truein 

{\small\leftskip 25truept\rightskip 25truept{\bf Abstract}\stdspace\theabstract

{\bf AMS Classification}\stdspace\theprimaryclass
\ifx\thesecondaryclass\relax\else; \thesecondaryclass\fi\par
{\bf Keywords}\stdspace \thekeywords\par}\vglue 7truept

}   %%%% end of definition of \makeagttitle

\ifplaintex
\font\phead=cmsl9 scaled 950
\font\pnum=cmbx10 scaled 913
\font\pfoot=cmsl9 scaled 950
\headline{\vbox to 0pt{\vskip -4.5mm\line{\small\phead\ifnum
\count0=\startpage ISSN 1472-2739 (on-line) 1472-2747 (printed)
\hfill {\pnum\folio}\else\ifodd\count0\def\\{ }% 
\ifx\theshorttitle\relax\thetitle\else\theshorttitle\fi\hfill{\pnum\folio}
\else\def\\{ and }{\pnum\folio}\hfill\ifx\theshortauthors\relax\theauthors
\else\theshortauthors\fi\fi\fi}\vss}}
\footline{\vbox to 0pt{\vglue 0mm\line{\small\pfoot\ifnum\count0=\startpage
\copyright\ \gtp\hfill\else
\agt, Volume \thevolumenumber\ (\thevolumeyear)\hfill\fi}\vss}}
\else
\font=cmsl9 scaled 1050
\font\lnum=cmbx10 
\font=cmsl9 scaled 1050
\makeatletter
\def\@oddhead{{\small\lhead\ifnum\count0=\startpage ISSN 1472-2739 
(on-line) 1472-2747 (printed)\hfill {\lnum\number\count0}\else\ifodd\count0
\def\\{ }\ifx\theshorttitle\relax \thetitle \else\theshorttitle\fi\hfill
{\lnum\number\count0}\else\def\\{ and }{\lnum\number\count0}
\hfill\ifx\theshortauthors\relax 
\theauthors\else\theshortauthors\fi\fi\fi}}\def\@evenhead{\@oddhead}
\def\@oddfoot{\small\lfoot\ifnum\count0=\startpage\copyright\ \gtp\hfill\else
\agt, Volume \thevolumenumber\ (\thevolumeyear)\hfill\fi}
\def\@evenfoot{\@oddfoot}
\makeatother
\fi
\let\maketitlepage\makeagttitle
\let\makeshorttitle\maketitlepage
\let\maketitle\maketitlepage

\newwrite\gtoutfile
\long\gdef\makeheadfile{  %%% start of definition of \makeheadfile
{\def\\{, }\def\s{ }
\immediate\openout\gtoutfile head.xxx
\immediate\write\gtoutfile{To: math@arxiv.org}
\immediate\write\gtoutfile{Subject: put OR rep NNNNN:ppppp}
\immediate\write\gtoutfile{--text follows this line--}
\immediate\write\gtoutfile{Proxy-for: \ifx\theasciiauthors\relax
\theauthors\else\theasciiauthors\fi\s<\ifx\theasciiemail\relax\theemail\else\theasciiemail\fi>}
\immediate\write\gtoutfile{\noexpand\\}
\immediate\write\gtoutfile{Authors: \ifx\theasciiauthors\relax
\theauthors\else\theasciiauthors\fi}
{\def\\{ }\immediate\write\gtoutfile{Title: \ifx\theasciititle\relax
\thetitle\else\theasciititle\fi}}
\immediate\write\gtoutfile{Subj-class: GT or SG, GR etc}
\immediate\write\gtoutfile{MSC-class: \theprimaryclass\ifx\thesecondaryclass\relax\else, \thesecondaryclass\fi}
\immediate\write\gtoutfile{Journal-ref: Algebr. Geom. Topol. \thevolumenumber\s
(\thevolumeyear) \startpage-\finishpage}
\immediate\write\gtoutfile{Comments: Published by Algebraic and
Geometric Topology at}
\immediate\write\gtoutfile{\s\s\s  http://www.maths.warwick.ac.uk/agt/AGTVol\thevolumenumber/agt-\thevolumenumber-\thepapernumber.abs.html}
\immediate\write\gtoutfile{\noexpand\\}
\immediate\write\gtoutfile{}
\ifx\theasciiabstract\relax
\immediate\write\gtoutfile{\theabstract}\else
\immediate\write\gtoutfile{\theasciiabstract}\fi
\immediate\write\gtoutfile{}
\immediate\write\gtoutfile{\noexpand\\}
\immediate\write\gtoutfile{}
\immediate\closeout\gtoutfile}}  %%% end of definition of \makeheadfile

\def\maketitlepage{\makeagttitle\makeheadfile}
\let\makeshorttitle\maketitlepage
\let\maketitle\maketitlepage

\def\ifplaintex{\expandafter\ifx\csname documentclass\endcsname\relax}

\def\gtp{{\mathsurround=0pt\it $\cal G\mskip-2mu$eometry \&\ 
$\cal T\!\!$opology $\cal P\!$ublications}}  % GT publications

\def\recd{{\small Received:\qua\receiveddate\ifx\reviseddate\relax
\else\qquad Revised:\qua\reviseddate\fi\par}} 

\def\lognumber#1{\def\thelognumber{#1}}
\def\volumenumber#1{\def\thevolumenumber{#1}}
\def\volumeyear#1{\def\thevolumeyear{#1}}
\def\papernumber#1{\def\thepapernumber{#1}}
\def\pagenumbers#1#2{\def\startpage{#1}\def\finishpage{#2}}
\def\published#1{\def\publishdate{#1}}

\def\received#1{\def\receiveddate{#1}}
\def\revised#1{\def\reviseddate{#1}}
\def\accepted#1{\def\accepteddate{#1}}

\long\def\asciiabstract#1{\long\def\theasciiabstract{#1}}
\def\asciikeywords#1{\def\theasciikeywords{#1}}

\let\\\par\let\thelognumber\relax\let\thevolumenumber\relax
\let\thepapernumber\relax\let\thevolumeyear\relax\let\startpage\relax
\let\finishpage\relax\let\publishdate\relax\let\receiveddate\relax
\let\reviseddate\relax\let\accepteddate\relax\let\theasciititle\relax
\let\theasciiauthors\relax
\let\theasciiabstract\relax\let\theasciikeywords\relax

\let\theasciiemail\relax

\ifplaintex
\font\logobig=cmssbx10 scaled 3836
\font\logomed=cmssbx10 scaled 2557
\else
\font\logobig=cmssbx10 scaled 4200
\font\logomed=cmssbx10 scaled 2800
\fi

\long\def\makeagttitle{   %%% start of definition of \makeagttitle
\count0=\startpage
\agt\hfill      %   Journal title (top left) 
\hbox to 45truept{\vbox to 0pt{\vglue -13truept{\logomed A\kern -.37em{\logobig 
T}\kern -.38em G}\vss}\hss}
\break
{\small Volume \thevolumenumber\ (\thevolumeyear)
\startpage--\finishpage\nl
Published: \publishdate}

\vglue .25truein

{\parskip=0pt\leftskip 0pt plus
1fil\def\\{\par\smallskip}{\Large\bf\thetitle}\par\medskip} \vglue
0.05truein

{\parskip=0pt\leftskip 0pt plus 1fil\def\\{\par}{\sc\theauthors}
\par\medskip}%
 
\vglue 0.03truein 

{\small\leftskip 25truept\rightskip 25truept{\bf Abstract}\stdspace\theabstract

{\bf AMS Classification}\stdspace\theprimaryclass
\ifx\thesecondaryclass\relax\else; \thesecondaryclass\fi\par
{\bf Keywords}\stdspace \thekeywords\par}\vglue 7truept

}   %%%% end of definition of \makeagttitle

\ifplaintex
\font\phead=cmsl9 scaled 950
\font\pnum=cmbx10 scaled 913
\font\pfoot=cmsl9 scaled 950
\headline{\vbox to 0pt{\vskip -4.5mm\line{\small\phead\ifnum
\count0=\startpage ISSN 1472-2739 (on-line) 1472-2747 (printed)
\hfill {\pnum\folio}\else\ifodd\count0\def\\{ }% 
\ifx\theshorttitle\relax\thetitle\else\theshorttitle\fi\hfill{\pnum\folio}
\else\def\\{ and }{\pnum\folio}\hfill\ifx\theshortauthors\relax\theauthors
\else\theshortauthors\fi\fi\fi}\vss}}
\footline{\vbox to 0pt{\vglue 0mm\line{\small\pfoot\ifnum\count0=\startpage
\copyright\ \gtp\hfill\else
\agt, Volume \thevolumenumber\ (\thevolumeyear)\hfill\fi}\vss}}
\else
\font=cmsl9 scaled 1050
\font\lnum=cmbx10 
\font=cmsl9 scaled 1050
\makeatletter
\def\@oddhead{{\small\lhead\ifnum\count0=\startpage ISSN 1472-2739 
(on-line) 1472-2747 (printed)\hfill {\lnum\number\count0}\else\ifodd\count0
\def\\{ }\ifx\theshorttitle\relax \thetitle \else\theshorttitle\fi\hfill
{\lnum\number\count0}\else\def\\{ and }{\lnum\number\count0}
\hfill\ifx\theshortauthors\relax 
\theauthors\else\theshortauthors\fi\fi\fi}}\def\@evenhead{\@oddhead}
\def\@oddfoot{\small\lfoot\ifnum\count0=\startpage\copyright\ \gtp\hfill\else
\agt, Volume \thevolumenumber\ (\thevolumeyear)\hfill\fi}
\def\@evenfoot{\@oddfoot}
\makeatother
\fi
\let\maketitlepage\makeagttitle
\let\makeshorttitle\maketitlepage
\let\maketitle\maketitlepage

\newwrite\gtoutfile
\long\gdef\makeheadfile{  %%% start of definition of \makeheadfile
{\def\\{, }\def\s{ }
\immediate\openout\gtoutfile head.xxx
\immediate\write\gtoutfile{To: math@arxiv.org}
\immediate\write\gtoutfile{Subject: put OR rep NNNNN:ppppp}
\immediate\write\gtoutfile{--text follows this line--}
\immediate\write\gtoutfile{Proxy-for: \ifx\theasciiauthors\relax
\theauthors\else\theasciiauthors\fi\s<\ifx\theasciiemail\relax\theemail\else\theasciiemail\fi>}
\immediate\write\gtoutfile{\noexpand\\}
\immediate\write\gtoutfile{Authors: \ifx\theasciiauthors\relax
\theauthors\else\theasciiauthors\fi}
{\def\\{ }\immediate\write\gtoutfile{Title: \ifx\theasciititle\relax
\thetitle\else\theasciititle\fi}}
\immediate\write\gtoutfile{Subj-class: GT or SG, GR etc}
\immediate\write\gtoutfile{MSC-class: \theprimaryclass\ifx\thesecondaryclass\relax\else, \thesecondaryclass\fi}
\immediate\write\gtoutfile{Journal-ref: Algebr. Geom. Topol. \thevolumenumber\s
(\thevolumeyear) \startpage-\finishpage}
\immediate\write\gtoutfile{Comments: Published by Algebraic and
Geometric Topology at}
\immediate\write\gtoutfile{\s\s\s  http://www.maths.warwick.ac.uk/agt/AGTVol\thevolumenumber/agt-\thevolumenumber-\thepapernumber.abs.html}
\immediate\write\gtoutfile{\noexpand\\}
\immediate\write\gtoutfile{}
\ifx\theasciiabstract\relax
\immediate\write\gtoutfile{\theabstract}\else
\immediate\write\gtoutfile{\theasciiabstract}\fi
\immediate\write\gtoutfile{}
\immediate\write\gtoutfile{\noexpand\\}
\immediate\write\gtoutfile{}
\immediate\closeout\gtoutfile}}  %%% end of definition of \makeheadfile

\def\maketitlepage{\makeagttitle\makeheadfile}
\let\makeshorttitle\maketitlepage
\let\maketitle\maketitlepage

\lognumber{25}
\volumenumber{1}
\volumeyear{2001}
\papernumber{25}
\pagenumbers{491}{502}
\received{26 October 2000}
\revised{7 May 2001}
\accepted{17 August 2001}
\published{10 September 2001}

\usepackage{amsmath,amscd,amssymb} %%,hyperref}
\usepackage[all]{xy}

\theoremstyle{plain}

\newtheorem{thm}[subsection]{Theorem}
\newtheorem{lem}[subsection]{Lemma}

\theoremstyle{remark}
\newtheorem*{rem}{Remark}

\theoremstyle{definition}

\numberwithin{equation}{section}

\newcommand{\cat}{\textup{cat}}
\newcommand{\Ext}{\textup{Ext}}

\begin{document}

\title{The product formula for Lusternik-Schnirelmann\\category}

\author{Joseph Roitberg}
\address{Department of Mathematics and Statistics\\Hunter College, 
CUNY\\695 Park Avenue, New York, NY 10021, USA}
%\address{Ph.D Program in Mathematics, The Graduate School and University
%Center, The City University
%of New York, 365 Fifth Avenue, New York, NY 10016-4309, USA}
\email{roitberg@math.hunter.cuny.edu}

\primaryclass{55M30}
\keywords{Phantom map, Mislin (localization) genus,
Lusternik-Schnirel\-mann category, Hopf invariant, Cuplength}
\asciikeywords{Phantom map, Mislin (localization) genus,
Lusternik-Schnirelmann category, Hopf invariant, Cuplength}

\begin{abstract}
If $C=C_\phi\,$ denotes the mapping cone of an essential phantom
map $\phi$ from the suspension of the Eilenberg-Mac Lane complex
$K=K(\mathbb{Z},5)\,$ to the $4$-sphere $S=S^4,$ we derive the
following properties: (1) The LS category of the product of $C$
with any $n$-sphere $S^n$ is equal to $3$; (2) The LS category of
the product of $C$ with itself is equal to $3$, hence is strictly
less than twice the LS category of $C$. These properties came to
light in the course of an unsuccessful attempt to find, for each
positive integer $m$, an example of a pair of $1$-connected
CW-complexes of finite type in the same Mislin (localization) genus
with LS categories $m$ and $2m.\,$If$\,\phi$ is such that its
$p$-localizations are inessential for all primes $p$, then by the
main result of \cite{R2}, the pair
$C_\ast=S\vee\Sigma^2K$,$\,C\,$provides such an example in the
case $m=1$.
\end{abstract}

\asciiabstract{If C=C_\phi denotes the mapping cone of an essential
phantom map \phi from the suspension of the Eilenberg-Mac Lane complex
K=K(Z,5) to the 4-sphere S=S^4 we derive the following properties: (1)
The LS category of the product of C with any n-sphere S^n is equal to
3; (2) The LS category of the product of C with itself is equal to 3,
hence is strictly less than twice the LS category of C. These
properties came to light in the course of an unsuccessful attempt to
find, for each positive integer m, an example of a pair of 1-connected
CW-complexes of finite type in the same Mislin (localization) genus
with LS categories m and 2m. If \phi is such that its p-localizations
are inessential for all primes p, then by the main result of
[J. Roitberg, The Lusternik-Schnirelmann category of certain infinite
CW-complexes, Topology 39 (2000), 95-101], the pair C_*, C where C_*=
S wedge \Sigma^2 K, provides such an example in the case m=1.}

\maketitle

\section{Introduction}

In this sequel to \cite{R2}, we record two additional curious
properties of the space $X$ described in the main Theorem of
\cite{R2}.

Recall that $X$ is the mapping cone of an essential map $\phi$
from $\Sigma K$, the suspension of the Eilenberg-Mac Lane
complex $K=K(\mathbb{Z},\,5)$, to$\,S=S^4$, the $4$-sphere.
Henceforth we use the notation $C=C_\phi$ for this space. Recall
from \cite{R2} that $\phi$ factors uniquely as
    $$\phi=\psi\circ r,$$
where
    $$r:\Sigma K\longrightarrow S^6_{(0)}$$
is a rationalization map and
    $$\psi: \,S^6_{(0)}\longrightarrow S.$$
\eject

Both $\phi$ and $\psi$ are phantom maps in the sense of \cite{Z}
(see also \cite[End of section 1]{R1}); they become inessential
when restricted to {\it finite} subcomplexes of their respective
domains. We study $D=D_\psi$ , the mapping cone of $\psi$, side
by side with $C$. The space $D$ is simpler (it is
finite-dimensional) and less interesting (it does not have finite
type) than $C$ and serves as a model for the latter. The two
results below apply to both spaces; the proofs in the (easier)
case of $D$ provide motivation and guidance for the proofs in the
case of $C$.

The first result states, in view of the main Theorem of
\cite{R2}, that the spaces $C,\,D\,$ satisfy the Ganea
``conjecture'' for all spheres. Thus:

%%Theorem 1.1
\begin{thm} For any $n>0,$
    $$\cat(C\times S^n)=3=\cat(D\times S^n).$$
\end{thm}

As an application of Theorem 1.1, we obtain a pair of spaces in the same Mislin
genus having LS categories $2$ and $3$. Namely, if $C_{\ast\,}$ (respectively
$D_\ast$) denotes the mapping cone of the trivial map with domain $\Sigma K$
(respectively $S_{(0)}^6$), and if $\phi$ is a {\it special} phantom map as in
\cite[Example 2]{R2} (that is the $p$-localization of $\phi$ is
inessential for every prime $p$), then the spaces $C\times S^n,\,C_\ast\times
S^n$ are in the same Mislin genus, as are the spaces $D\times
S^n,\,D_\ast\times S^{n\text{ }}$, and
    $$\cat(C\times S^n)=3, \qquad \cat(C_\ast\times S^n)=2;$$
    $$ \cat(D\times S^n)=3, \qquad \cat(D_\ast\times S^n)=2.$$

This application is only moderately interesting since we may
generalize Example 2 of \cite{R2} to get an example of a pair of
spaces in the same Mislin genus having LS categories $m$ and
$m+1$ for {\it any} $m>0$. Indeed, we need only replace $K$ by
$K(\mathbb{Z},4m+1)\,$and $S$ by $\mathbb{H}P^m$, quaternionic
projective $m$-space, and argue as in  \cite{R2}, applying
appropriate results of Iwase \cite{I2}; details are omitted. What
would be more interesting is to find a pair of spaces in the same
Mislin genus with LS categories $m$ and $2m$ {\it for any} $m>0$. For
according to a result of Cornea \cite{C}, with improvement by
Stanley \cite{S1} (see also F\'elix-Halperin-Thomas \cite{FHT}),
    \begin{equation*} \cat(X)\leq 2{\cdot}\sup_{p} \bigl(\cat(X_{(p)})\bigr),
    \tag{1.1} \end{equation*}
where $X$ is a $1$-connected CW-complex of finite type and $X_{(p)}$ is its
$p$-localiza\-tion.

%%3

Examples of the type just mentioned would show that the inequality (1.1) is
sharp for all $m>0$. To that end, our original thought was to try the
$m$-fold \ product spaces \ $C^m_\ast,C^m\,$(or $D^m_\ast,D^m$). However, this
attempt fails, as the next result shows.

%%Theorem 2
\begin{thm}
$\cat(C_{\,}\times C)=3=\cat(D\times\,D).$
\end{thm}

Thus $C$ is a $1$-connected CW-complex of finite type such that
$\cat(C\times C)$ is strictly smaller than $2{\cdot}\cat(C)$\,.
Examples of strict inequality in the product formula $\cat(A\times
B)\leq\cat(A)+\cat(B)$ go back to Fox \cite{F}. However, the
two factors in Fox's example are distinct spaces: they are Moore
spaces $S^2\cup_p\,e^3,S^2\cup_q\,e^3$ with respect to distinct
primes $p$ and $q$. See also Ganea-Hilton \cite{GH} for a
generalization of Fox's example as well as the observation that
$\cat(K(\mathbb{Q},1)\times K(\mathbb{Q},1))\leq 3$\, (which is
strictly smaller than $2{\cdot}\cat(K(\mathbb{Q},1))=4$);  of
course, $K(\mathbb{Q},1)$, unlike $C$, is neither $1$-connected
nor of finite type. An alternate derivation of the inequality
$\cat(K(\mathbb{Q},1)\times K(\mathbb{Q},1))\leq 3$ is possible
along the lines of the first part of the proof below of Theorem
1.2.

A more recent example of strict inequality in the product formula, similar
in flavor to Fox's example (in both situations, $A\vee B$ is homotopy
equivalent to $A\times B$), is contained in F\'{e}lix-Halperin-Lemaire
\cite{FHL}.  In this example, $A$ is a $1$-connected CW-complex whose
(reduced) homology groups are all finite and $B$ is a $1$-connected,
rational CW-complex. The main result of \cite{FHL} states that --- unlike
the situation when $A$, $B$ are both $K(\mathbb{Q},1)$ --- $\cat(A\times
B)=\cat(A)+\cat(B)$ whenever $A,B$ are rationalizations of
\emph{1-connected} CW-complexes of finite type.

The landscape changed radically when Iwase \cite{I1} discovered
examples of strict inequality in the product formula in which one
of the factors (but not both!) is a sphere. Systematic approaches
to Ganea's ``conjecture'', with many more examples, have been
developed since then; see \cite{I2}, \cite{S1},\ Vandembroucq
\cite{V},\ Stanley \cite{S2} and\ Harper \cite{H}.

We remark that the examples in  \cite{H} and \cite{S2}, \
illustrating the strict inequality  $\cat(A\times A)<2{\cdot}\cat(A)$
for $A$ $1$-connected of finite type are quite different from
those in Theorem 1.2. In each of Harper's examples, $A$ is the
mapping cone of a suitable map from a sphere to another sphere and
$\cat(A)=2$, $\cat(A\times A)\leq 3$. In Stanley's example, $A$ is
the mapping cone of a suitable map from a sphere to a bouquet of
two spheres and $\cat(A)=2=\cat(A\times A)$. In contrast to the
situation in Theorem 1.1, the examples in \cite{H} and \cite{S2}
do \emph{not} satisfy the Ganea ``conjecture''. In fact, $\cat(A\times
S^n)=2$ for all $n\geq 2$ in two of the three examples in
\cite{H} and  $\cat(A\times S^n)=2$ for all $n\ge1$ in the third
example in \cite{H} and in the example in \cite{S2}.

The proofs of Theorems 1.1 and 1.2 are contained in the next
section. The proof of Theorem 1.1 builds on the proof of the main
Theorem of \cite{R2}, bringing in more of the Hopf invariant
technology of \cite{I2}.

We also take the opportunity to correct a mis-statement in the last
paragraph on p. 99 of \cite{R2}; I thank Jianzhong Pan for bringing
the error to my attention.

As for the proof of Theorem 1.2, the
inequalities
    $$ \cat(C\times C)\leq\,3, \qquad \cat(D\times D)\leq\,3$$
follow from classical obstruction theory arguments. Similar, but slightly
more elaborate, arguments lead to the inequalities
    $$ \cat(C^k)\leq k+1, \qquad \cat(D^k)\leq k+1$$
for $3\leq k\leq 6$ as well, but
these arguments break down for $k=7$. To establish the opposite
inequalities
    $$ \cat(C\times\,C)\geq\,3, \qquad \cat(D\times D)\geq\,3,\,$$
we exploit techniques of \cite{I2} together with
cuplength arguments. In fact, we show that $D\times D\,$ has
cuplength $3$. This latter fact fails for $C\times C$ (whose
cuplength is $2$). To remedy this failure, we introduce a notion
of \textit{refined cup products} - which in the case of $C\times
C$ consists of a blend of ordinary cohomology and cohomotopy - and
show that $C\times C\,$has non-zero refined cup products of length
$3$, thereby leading to the desired inequality.

The latter approach also leads to an alternate, albeit less
direct, proof of Theorem 1.1, as well as of the main Theorem of
\cite{R2}.

%%Section 2: proofs
\section{Proofs}

Before giving the proof of Theorem 1.1, which closely follows the
lines of the proof of the main Theorem of \cite{R2}, we wish to
point out that throughout the last paragraph on page 99 of
\cite{R2}, the space $S^6_{(0)}\,$ should be replaced by $\Sigma
K$. To see that the argument beginning in that paragraph and
concluding on page 100 is valid for $\Sigma K$, it is only
necessary to observe that we may identify the map
    $$[\Sigma K,S^7]\longrightarrow[\Sigma K,E^2(\Omega S^4)]$$
induced by the inclusion of the bottom $S^7$ into $E^2(\Omega S^4)$ with
the analogous map
    $$[S^6_{(0)},S^7]\longrightarrow[S^6_{(0)},E^2(\Omega S^4)],$$
denoted in \cite{R2} by $\beta_1$; here $E^2(\Omega S^4)$ is a space
homotopy equivalent to the 2-fold join $\Omega S^4 \ast \Omega
S^4$. Indeed, we have a commutative diagram
    \[
    \begin{CD}
    [S^6_{(0)},S^7] @>>> [S^6_{(0)},E^2(\Omega S^4)] \\
    @VVV @VVV\\
    [\Sigma K,S^7] @>>> [\Sigma K,E^2(\Omega S^4)]
    \end{CD} \]
where the vertical maps are induced by rationalization $r:\Sigma K
\to\ S^6_{(0)}\,$ and are isomorphisms; see \cite[(2.5)]{R2}.

We now give the proof of Theorem 1.1.

\begin{proof}[Proof of Theorem 1.1]
It is clear that
    $$ \cat(A\times S^m)\leq 3, \qquad A=C \textup{ or } D,$$
for any $m\geq 1$. To prove that
    $$ \cat(A\times S^m)=3,$$
we appeal to a special case of \cite[Theorem 3.9(2)]{I2} according to
which it suffices to establish the following variant of (iii) on page 98
of \cite{R2}:
    \begin{equation*} (E^2(\Omega i)\ast j_{m-1})\circ\Sigma^mH_1(\alpha)
    \tag{iii$'$} \end{equation*}
is essential, where the left-most map is the join of the map induced
by the inclusion $i:S\subset A\,$ with the canonical map
$j_{m-1}:S^{m-1}\to\ \Omega S^m$ , $\alpha=\phi$ ($\,$respectively
$\psi$) if $A=C$ (respectively $D$), and $H_1(\alpha)$ is the
Berstein-Hilton-Hopf invariant of $\alpha$ as generalized by
Iwase (\cite[Definition 2.4]{I2}) and Stanley
(\cite[Definition 3.4]{S1}).

Now the argument on page 99 of \cite{R2} prior to the last
paragraph actually shows that $\Sigma^mH_1(\alpha)\,$ is essential
for any $m \geq 0$. Moreover the subsequent argument (suitably
amended; see above) shows that $E^2(\Omega i)\,$ has a left
homotopy inverse. Since $\Sigma j_{m-1}\,$ has a left homotopy
inverse$,\,$it follows that $E^2(\Omega i)\ast j_{m-1}$ too has a
left homotopy inverse, and (iii$'$) is verified.

In the case $A=D$, it is even permissible to ``ignore'' the map $E^2(\Omega
i)\ast j_{m-1}$; indeed, that map induces a monomorphism from $[S^6_{(0)}\ast
S^{m-1},E^2(\Omega S^4)\ast S^{m-1}]$ to $[S^6_{(0)}\ast S^{m-1},E^2(\Omega
D)\ast\Omega S^m]$ since $S^6_{(0)}\,$has dimension $7$. See also
\cite[Remark 3.10]{I2}.
\end{proof}

\rem
Since any rationalization map $r:\Sigma K\rightarrow
S^6_{(0)}$ is a suspension, it follows from \cite[Proposition 2.11(1)]{I2}
that
    $$H_1(\varphi)=H_1(\psi)\circ r.$$

We next give the proof of Theorem 1.2.

\begin{proof}[Proof of Theorem 1.2]
To prove
    \begin{equation*}  \cat(B)\leq 3, \qquad
    B=C\times C \textup{ or }  D\times D, \tag{2.1} \end{equation*}
it suffices to show that the Ganea fibration
    \begin{equation*}  P^3(\Omega B)\longrightarrow B, \tag{2.2} \end{equation*}
with fiber
    $$E^4(\Omega B)\simeq\Omega B\ast\Omega B\ast\Omega B\ast\Omega
    B,$$
admits a section; see \cite[Theorem 2.1]{I2}. For either
value of $B$ in (2.1), we easily verify that $E^4(\Omega B)$ is
14-connected - more precisely, $E^4(\Omega B)$ has the homotopy
type of a CW-complex of the form
    $$(S^{15}\vee \cdots \vee S^{15})\cup\bigl\{\hbox{cells
    of dimension} \geq 18\bigr\}.$$
Next observe that
    $$H_n(D\times D)=0$$
for $n>14$; note that unless the coefficients are
explicitly indicated, homology and
cohomology groups are understood to be integral. Thus
    $$H^n(D\times D;\pi_{n-1}(E^4(\Omega(D\times D)))=0$$
for all $n$. Since the obstructions to a section of (2.2) in the case
$B=D\times D$ lie in these cohomology groups, it is clear that a section
exists.

Two remarks:

\begin{enumerate}
\item An alternate approach to the existence of a section of (2.2) in this
case, but not in the case $B=C\times C$, is via a slight variant
of \cite[Theorem 5.5]{I2} --- details are omitted.
\item  We do not assert that there is a unique section of (2.2)
in the case $B=D\times D$. In fact, the (vertical homotopy classes
of) sections are classified by
    $$H^{15}(D\times D;\pi_{15}(E^4(\Omega(D\times D)))\approx \Ext
    (\mathbb{Q},\mathbb{Z}\oplus \cdots \oplus\mathbb{Z}),$$
which has the cardinality of the continuum.
\end{enumerate}

To deal with the case $B=C\times C$, first note that a homology
decomposition of $K,$ say
    $$S^5=K[5]\subset \cdots \subset K[m]\subset \cdots \subset K,$$
where each $K[m]\,$ is a finite
complex, naturally induces a homology decomposition of $C,$
    $$S=C[4]\subset \cdots \subset C[m]\subset \cdots \subset C,$$
where each $C[m]\,$ is likewise a finite complex. Explicitly, we take
$C[m]$ to be the mapping cone of the restriction of $\phi$ to $\Sigma
K[m-1]$. Thus $C\times C$ is the ascending union of the finite
subcomplexes $C[m]\times C[m]$, about which we make two claims:

%% second internal enumeration
\begin{enumerate}
\item  $\cat(C[m]\times C[m])=2$, $m\geq 4$. In fact, the homotopy equivalence
$C[m]\simeq S\vee\Sigma^2K[m-1]$, $m\geq 4$, resulting from the fact that
$\phi$ is phantom, implies that $\cat(C[m])=1$, hence that $\cat(C[m]\times
C[m])=2$, $m\geq 4$.
\item For each $m$, $H_n(C[m]\times C[m])$ is finite for all $n>14$.
This follows immediately from the fact that $H_n(C\times C)$ is
finite for all $n>14$.
\end{enumerate}
%% end second internal enumeration
%%7
It follows from (1) that for each $m$, the inclusion $C[m]\times
C[m]\subset C\times C$ lifts to the total space
$P^3(\Omega(C\times C))$ of the fibration (2.2). Moreover, the
number of such lifts is finite since, thanks to (2), the
cohomology groups $H^n(C[m]\times C[m];\pi_n(E^4(\Omega(C\times
C)))$ are finite for all $n$. A variant of the classical argument
utilized in the first paragraph of page 100 of \cite{R2} then
leads to the conclusion that (2.2) admits a section in the case
$B=C\times C$.

We now verify that
    $$ \cat(B)\geq 3, \qquad B=C\times C \textup{ or } D\times D,$$
beginning with the
case $B=D\times D$. Our strategy is to detect a non-zero $3$-fold
cup product in the integral cohomology ring of $D\times D.\,$More
precisely, if $u$ and $v$ are the canonical generators (up to
sign) of $H^4(D\times D)$, we show that
    \begin{equation*} u^2v\neq 0. \tag{2.3} \end{equation*}

We first prove the weaker statement
    \begin{equation*}   u^2\neq 0 \tag{2.4} \end{equation*}
and then generalize the proof of (2.4) to
obtain (2.3). Let $u_\ast$ be the restriction of $u$ to
$D\times\ast\subset D\times D$. To prove (2.4), it suffices to
prove
    \begin{equation*}  u_\ast^2\neq 0. \tag{2.5} \end{equation*}
Applying \cite[Theorem 4.1, Remark 4.3]{I2} to the cofibration
sequence
    $$S^6_{(0)}\longrightarrow S\longrightarrow D,$$
we infer that the reduced diagonal map
    $$\Delta'=\Delta'_D:D\longrightarrow D\wedge D$$
factors as the composite
    \begin{equation*} (i\wedge i)\circ(e\wedge e)\circ\Sigma
    H_1(\psi)\circ q:D \longrightarrow\ D\wedge D\ , \tag{2.6} \end{equation*}
where
\begin{enumerate}
\item $q: D\to\ D/S=S^7_{(0)}$ is the collapsing map,
\item $\Sigma H_1(\psi):S^7_{(0)} \to\ \Sigma E^2(\Omega S)\simeq\Sigma\Omega
S\wedge\Sigma\Omega S $ is the suspension of the
Berstein-Hilton-Hopf invariant of $\psi$, and
\item $(i\wedge i)\circ(e\wedge e)$ is induced by
the evaluation map $e:\Sigma\Omega S \rightarrow S$  and the
inclusion map $i: S\subset D$.
\end{enumerate}

Abbreviating $K(m)=K(\mathbb{Z},m)\,$ and viewing $u_\ast$ as (the
homotopy class of) a map
    $$u_\ast:D\longrightarrow K(4),$$
we may view $u^2_\ast$ as the composite
    \begin{equation*} \mu\circ(u_\ast\wedge u_\ast)\circ\Delta':D\longrightarrow
    D\wedge D \longrightarrow  K(4)\wedge K(4) \longrightarrow K(8),
    \tag{2.7} \end{equation*}
where $\mu$ represents a generator of $H^8(K(4)\wedge K(4))$. From (2.6)
and (2.7), we see that $u_\ast^2$ may be represented as the composite
    \begin{equation*}  \mu\circ(j\wedge j)\circ(e\wedge e)\circ\Sigma
    H_1(\psi)\circ q:D\longrightarrow K(8), \tag{2.8} \end{equation*}
where $j \wedge j$ is induced by a generator $j$ of $H^4(S)$ while the
other maps in the composite have been defined above.

Reasoning as in the proof of \cite[(2.3)]{R2} shows that $(e\wedge
e)\circ\Sigma H_1(\psi)$ is non-zero. It is then routine to
argue that the composite (2.8) is itself non-zero. Thus (2.5),
and with it (2.4), is proved.
%%$\,$ \

Next let $u_{\bullet},\,v_{\bullet}$ be the images of $u,v$ under
the cohomology map induced by the inclusion $D\times S\subset
D\times D$. Of course, the restriction of \ $u_{\bullet}$ to
$D\times\ast\subset D\times S\,$ is just $u_\ast$.
To prove (2.3), it suffices to prove
    \begin{equation*}  u_{\bullet}^2 v_{\bullet} \neq 0. \tag{2.9} \end{equation*}

To that end, we apply \cite[Corollary 4.1.1]{I2} to the cofibration
sequence
    $$S^{10}_{(0)}\longrightarrow (D\times\ast)\cup(S\times S)\to\
    D\times S$$
to infer that the reduced diagonal map
    $$\Delta'\wedge 1_S=q'\circ(\Delta'\times 1_S):D\times S\to\
    (D\wedge D)\times S\longrightarrow  D\wedge D\wedge S,$$
where $q'\,$is the collapsing map, factors as
    \begin{equation*} (i\wedge i\wedge 1_S)\circ(e\wedge e\wedge
    1_S)\circ\Sigma^5H_1(\psi)\circ q^{!}: D\times S \longrightarrow
    D\wedge D\wedge S \tag{2.10} \end{equation*}
where
\begin{enumerate}
\item $q^{!}: D\times S\to\ D\times
S/(D\times\ast)\cup(S\times S)=(D/S)\wedge S$ is the collapsing
map,
\item $D\times S/(D\times\ast)\cup(S\times S)=(D/S)\wedge S = S^{11}_{(0)}$,
\item $\Sigma^5H_1(\psi): S^{11}_{(0)} \to\ \Sigma\Omega S\wedge\Sigma\Omega
S\wedge S $
is the 5-fold iterated suspension of $H_1(\psi)$,
\item $e\wedge e\wedge 1_S: \Sigma\Omega S\wedge\Sigma\Omega S\wedge S
\to\  S\wedge S\wedge S$ is induced by $e$, and
\item $i\wedge i \wedge 1_S:S\wedge S\wedge S \to\ D\wedge D\wedge
S$ is induced by $i$.
\end{enumerate}

From (2.10), we see that $u^2_{\bullet} v_{\bullet}$ may be
represented as the composite
    \begin{equation*}  \nu\circ(j\wedge j\wedge j)\circ(e\wedge e\wedge
    1_S)\circ\Sigma^5H_1(\psi)\circ q^{!}:D\times S \longrightarrow K(12)
    \tag{2.11} \end{equation*}
where $\nu$ represents a generator of $H^{12}(K(4)\wedge K(4)\wedge K(4))$.

Arguing as with (2.8), we readily verify that the composite
(2.11) is non-zero. Thus (2.9), and with it (2.3), is
proved.

Finally, we adapt the preceding argument to prove the inequality
    $$ \cat(C\times C)\geq 3.$$
If $u,v$ now stand for the canonical generators (up to sign) of
$H^4(C\times C)$, the analogs of (2.3) and (2.4) both fail. In
fact, as the cohomology of $\,C$ is indistinguishable from that of
$C_\ast$, the cuplength of $C$ is precisely $1$ and the cuplength
of $C\times C$ is precisely $2$. One might try replacing
cohomology by cohomotopy since $\pi^8(C)=[C,S^8]\,$(as also
$\pi^8(D)=[D,S^8]\approx H^8(D))\,$ turns out to be
set-theoretically equivalent to $\Ext(\mathbb{Q},\mathbb{Z})$ as
noted in Lemma 2.1 below. Unfortunately, the Puppe sequence
    \begin{equation*} [\Sigma
    K,S]\longleftarrow[S,S]\longleftarrow[C,S]\longleftarrow[\Sigma^2K,S]
    \tag{2.12} \end{equation*}
associated to the cofibration sequence $\Sigma K\rightarrow
S\rightarrow C$ gives that
    $$\pi^4(C)=[C,S]=0$$
since the left-most map in (2.12) sends
$[S,S]$ monomorphically to the (infinite cyclic) group generated
by $\phi$, and $[\Sigma^2K,S]=0$ by \cite[Theorem D]{Z} (see also
\cite[Theorem 4.2]{R1}). In other words, the cohomology classes
$u,v$ are not compressible into $S$.

While neither cohomology alone nor cohomotopy alone supports a
cuplength argument necessary to establish the inequality
$\cat(C\times C)\geq 3$, the following compromise between
cohomology and cohomotopy provides a way out of the dilemma. Fix a
CW-model for $K(m)$ such that the $m$-skeleton is $S^m$ and the
$N$-skeleton $K(m)_N$ has only finitely many cells for all $N\geq
m$.

%%10

\begin{lem}
\begin{enumerate}

\item For $N\geq 9$, the inclusions $K(4)_N\subset
K(4)_{N+1}\subset K(4)$ induce a commutative diagram
    \begin{equation*}
    \begin{CD}
    [S,K(4)_N] @<<< [C,K(4)_N] \\
    @VVV @VVV \\
    [S,K(4)_{N+1}] @<<< [C,K(4)_{N+1}] \\
    @VVV @VVV \\
    [S,K(4)] @<<< [C,K(4)]
    \end{CD} \tag{2.13} \end{equation*}
with each of the sets in (2.13) equivalent to $\mathbb{Z}$, and
each of the arrows a bijection;

\item For $N\geq 8$, the inclusion $K(8)_N\subset
K(8)_{N+1}$ induces a commutative diagram
    \begin{equation*} \begin{CD}
    [C,K(8)_N]\hspace{0.4cm} @<<< [\Sigma^2K,K(8)_N] \\
    @VVV @VVV \\
    [C,K(8)_{N+1}]  @<<< [\Sigma^2K,K(8)_{N+1}]
    \end{CD} \tag{2.14} \end{equation*}
with each of the sets in (2.14) equivalent to $\Ext(\mathbb{Q},\mathbb{Z})$,
and each of the arrows a bijection;

\item For $N\geq 12$, the inclusion $K(12)_N\subset
K(12)_{N+1}$ induces a commutative diagram
    \begin{equation*} \begin{CD}
    [C\times S,K(12)_N] @<<< [\Sigma^6K,K(12)_N] \\
    @VVV @VVV \\
    [C\times S,K(12)_{N+1}] @<<< [\Sigma^6K,K(12)_{N+1}]
    \end{CD} \tag{2.15} \end{equation*}
with each of the sets in (2.15) equivalent to
$\Ext(\mathbb{Q},\mathbb{Z})$, and each of the arrows a bijection.
\end{enumerate}
\end{lem}

\begin{proof}
To prove the first part of (1), consider the commutative diagram
with exact rows
    \begin{equation*}
    \minCDarrowwidth1.5pc
%%  \begin{small}
    \begin{CD}
    [\Sigma K,K(4)_N] @<<< [S,K(4)_N] @<<< [C,K(4)_N] @<<<
    [\Sigma^2K,K(4)_{N}] \\
    @VVV @VVV @VVV @VVV \\
    [\Sigma K,K(4)_{N+1}] @<<< [S,K(4)_{N+1}] @<<< [C,K(4)_{N+1}]
    @<<< [\Sigma^2K,K(4)_{N+1}]
    \end{CD}
%%  \end{small}
    \tag{2.16} \end{equation*}
induced by the cofibration sequence $\Sigma K\rightarrow
S\rightarrow C$. Notice that if $M\geq 9,$
    $$\pi_{_7}(K(4)_M)=0=\pi_{_8}(K(4)_M).$$
Therefore, by \cite[Theorem D]{Z} (see also \cite[Theorem 4.2]{R1}),
    $$[\Sigma K,K(4)_M]=0=[\Sigma^2K,K(4)_M]$$
and (2.16) collapses to the upper square in (2.13). The proof of the rest
of (1), and also (2) and (3), are similar. Observe that the diagram in (3)
is induced by the cofibration sequence
    \begin{equation*} \Sigma^5 K \longrightarrow (C\times\ast)\cup(S\times
    S)\longrightarrow C\times S. \tag*{\qed} \end{equation*}
\def\qed{}\end{proof}

Define $u_{\ast N},\,u_{\bullet N},v_{\bullet N}$ to be the unique
compressions of $u_\ast,u_{\bullet},v_{\bullet}$ into
$K(4)_N$, $N\geq 9$, guaranteed by part (1) of Lemma 2.1. Further,
define two maps
    $$ \mu_N:K(4)_N\wedge K(4)_N \longrightarrow K(8)_{2N+1},$$
and
    $$\nu_N:K(4)_N\wedge K(4)_N\wedge K(4)_N \longrightarrow
    K(12)_{3N+1}$$
as the maps uniquely determined by $\mu$, $\nu$. It may be assumed
that these maps are cellular. We define the \textit{refined cup
square} $u^2_{\ast N}$ as the element represented by the
composite
    $$\mu_N\circ(u_{\ast N}\wedge u_{\ast
    N})\circ\Delta'_C:C\longrightarrow
    C\wedge C\longrightarrow K(4)_N\wedge K(4)_N \longrightarrow
    K(8)_{2N+1}.$$
We define the \textit{refined 3-fold cup product}, $u_{\bullet
N}^2v_{\bullet N}$ similarly, using $\nu_N$. Observe that
$u_{\ast N}^2$ compresses uniquely into $S^8$ by part (2) of
Lemma 2.1. In other words, $u_{\ast N}\,$is essentially a
$4-$dimensional cohomology class while $u^2_{\ast N}$ is
essentially an $8-$dimensional cohomotopy class. Similarly, by
part (3) of Lemma 2.1, $u^2_{\bullet N}v_{\bullet N}$ is
essentially a $12-$dimensional cohomotopy class.

We now invoke the analog of (2.8) with $D$ replaced by $C$,
$S^7_{(0)}\,$replaced by $\Sigma^2K$ and $K(m)$ replaced by
$K(m)_N,N\geq 9,$to conclude, with the help of Lemma 2.1, that
    $$u^2_{\ast N}\neq 0,$$
that is the \textit{refined  cuplength} of $C$
is at least $2$. Moreover, arguing with the appropriate analog of
(2.11), we find, once again with the help of Lemma 2.1, that
    $$u^2_{\bullet N}v_{\bullet N}\neq 0,$$
that is the refined cuplength of $C\times C$ is at least $3$.

The classical argument showing that cuplength is a lower bound for
cat readily generalizes to show that refined cuplength is also a
lower bound for cat. Hence $\cat(C\times C)$ is at least $3$,
completing the proof of Theorem 1.2.
\end{proof}

This research was supported in part by a grant from the City
University of New York  PSC-CUNY Research Award Program.

\bibliographystyle{amsalpha}

\Addresses

\recd

\end{document}